\documentclass[a4paper, 12pt]{article}
\usepackage{amsmath}
\usepackage{amsfonts}
\usepackage{amsthm}
\usepackage{mathrsfs}
\usepackage[T1]{fontenc}
\usepackage{amsmath}
\usepackage{amsfonts}
\usepackage{amsthm}
\usepackage{mathrsfs}
\usepackage{amssymb}
\usepackage{pst-node}
\usepackage{tikz-cd}

\newtheorem{prop}[]{Proposition}
\newtheorem{lem}[]{Lemma}

\newtheorem{cor}{Corollary}

\newtheorem{rem}[]{Remark}

\newcommand{\Rset}{\mathbb{R}}

\begin{document}

\title{Background driving distribution functions and series representations for
log-gamma selfdecomposable random variables}
\author{ Zbigniew J. Jurek (University of Wroc\l aw)}
\date{October 10 , 2021;  \ \  to appear in \emph{Theor. Probab. Appl.} 2022.}
\maketitle

\begin{quote} \textbf{Abstract.}
We identify the background driving  distribution functions (BDDF) for selfdecomposable distributions (random variables). For log-gamma variables and their background driving variables we find their series representations. An innovation variable for Bessel-K distribution is given as a compound Poisson variable.
\end{quote}


Primary 60E07, 60E10, 60H05; Secondary  60G51, 60B10.


\medskip
\emph{Abbreviated title: Background driving  distribution functions}

\medskip
\medskip
\medskip
Author's address:  Institute of Mathematics , University of Wroc\l aw,
Pl. Grunwaldzki 2/4 ,  50-384 Wroc\l aw ,  Poland. \\ www.math.uni.wroc.pl/$\sim$zjjurek ;
e-mail: zjjurek@math.uni.wroc.pl

\newpage

\textbf{0. INTRODUCTION.}

The L\'evy class, $L$, of \emph{selfdecomposable probability distributions (or random variables or characteristic functions)} is obtained as the limit of a linear normalization of partial sums of sequences of independent random variables. It constitutes  a proper subclass of the class, ID, of all infinitely divisible distributions; see the books by Gnedenko and Kolmogorov (1954), Sections 29-30, Feller (1966), Chapter XVII or Lo\'eve (1963), Section 23.  More recently in Bradley and Jurek (2016) it was proved that $L$ can be obtained by assuming only that we have strongly mixing sequences, keeping in mind that sequences of independent random variables are strongly mixing ones.

$L$ is a very large class and includes among others:  stable,  gamma, log-gamma, log-normal, $\chi^2$, Student-t, $\log|t|$, Snedecor F, hyperbolic variables and many other probability measures; cf. for instance Jurek (1997) and references therein. Moreover, the Thorin class, EGGC, of extended generalized gamma convolutions, is a subset of the class $L$; \ cf. Bondesson (1992), Chapter 7.

\medskip
The class $L$ has the following two equivalent characterizations. Namely, a selfdecomposable random variable $X$   admits:

 (a) \emph{a random integral representation}
\begin{equation}
X:= \int_0^\infty e^{-s}dY(s) \  \ \mbox{with} \   \ \mathbb{E}[1+\log|Y(1)|]<\infty,
\end{equation}
with respect to some (uniquely determined) L\'evy process $Y(t), t\ge 0$; and,

(b) \emph{a decomposition property}
\begin{equation}
\forall ( 0<c<1) \ \exists \ ( X_c  \bot   X )\  \ X  \stackrel{d}{=}\,cX+X_c,
\end{equation}
where $\bot$ means independence, and $X_c$ is referred to as an \emph{innovation} in the context of autoregressive sequences. Cf. Jurek and Vervaat (1983) or Jurek and Mason (1993), Chapter 3.

\medskip
\textbf{Terminology.} We refer to the process  $(Y(t), t\ge 0)$ in (1)  as \emph{the background driving L\'evy  process} (BDLP) of $X$, to $Y(1)$ -- as \emph{the background driving random variable} (BDRV) of $X$, and to the probability distribution function  $G_X(a):= P(Y(1)\le a), a \in\Rset$, as \emph{the background driving distribution function} (BDDF) of  $X$. Finally, the characteristic function $\psi$ of the random variable $Y(1)$ is called \emph{the background driving characteristic function} (BDCF) of the characteristic function $\phi$ of $X$.

\medskip
\begin{rem}
\emph{ a)  In  Jurek and Vervaat (1983), pp. 252-253,  using the remainders $\{X_c:  0\le c\le 1\}$ from (2), a cadlag process $Z(t), 0\le t< \infty,$ with independent increments was constructed  such that  $Y(t):=\int_{(0,t]}e^sdZ(s)$ is the BDLP of  $X$ from (1).}

 \emph{ b) In  Jeanblanc, Pitman and Yor (2002),  Theorem 1,  identified the  BDLP of X as
$
Y(t):= \int_{e^{-t}}^1\, r^{-1}dV(r), \ \   t\ge 0 , \ \   V(1):=X, \  (\mbox{X is from  (1))}$,
for an additive  1-self-similar process  $V(r), r\ge0$. In other words, selfdecomposable $X$ can be inserted into an additive 1-similar process. In fact,  we could have any H-self-similar process. [One should  keep in mind that the process V(r) from above is denoted as  $X_r$ in  Jeanblanc, Pitman and Yor (2002).]
}
\end{rem}

For  the identification  of BDDFs  $G_X(a):=P(Y(1)\le a)$, in principle,  we may use   $Y(1)$   from Remark 1  a) or b).
However, in this note  we will  utilize the characteristic functions $\phi$ of $X$  and $\psi$ of $Y(1)$, respectively, cf. Proposition 1, below. The  proposed method of finding BDDF is  illustrated by some explicit  examples; (cf. Lemma 1 and Section III).

Finally,  for log-gamma  and Bessel-K variables, the corresponding BDDFs  and innovation variables  are  described as  infinite series of some exponential variables (Propositions 2 - 4 and Corollary 5. These potentially  may help in numerical simulations of selfdecomposable variables ; cf. discussions in Hosseini (2019) and some references therein. In Corollary 4 there is a formula for the error function.

\medskip
\textbf{I. THE  RESULTS.}

\textbf{I. 1.} \emph{Background driving distribution functions (BDDF).}

Here is the BDDF formula for the selfdecomposable variable $X$ with the characteristic function $\phi$ :
\begin{prop} Let  a selfdecomposable $X$ have the integral representation (1) and $\phi(t):= \mathbb{E}[e^{itX}]$ be its characteristic function. Then $\phi$  is differentiable  ($t\neq 0$) and  the function
\begin{equation}
G_{X}(a):=P(Y(1)\le a)=\frac{1}{2}-\frac{1}{\pi}\int_0^\infty \Im(\exp[-ita +  t \frac{\phi^\prime(t)}{\phi(t)}])\,\frac{dt}{t},  \ \  a \in C_{G_{X}};
\end{equation}
is  the  BDDF  of $X$. [Above $C_{G_{X}}$ denotes the continuity points of  the function $G_{X}$ and $\Im z $ stands for imaginary part of $z$.]
\end{prop}
In the case of symmetric selfdecomposable variables we get the following:
\begin{cor}
Suppose that   a selfdecomposable  random variable X is a symmetric random variable, then its BDDF is given by
\begin{equation}
G_{X}(a)=\frac{1}{2}+\frac{1}{\pi}\int_0^\infty \exp[ t \frac{\phi^\prime(t)}{\phi(t)}]\, \sin(ta)\,\frac{dt}{t}, \ \ \ a \in C_{G_X}.
\end{equation}
\end{cor}
 Using  b)  from Remark 1 and the relation between the characteristic function of a selfdecomposable variable and its background driving characteristic function ( see (10) below) we have
\begin{equation}
\mathbb{E}[\exp(it\,\int_{e^{-1}}^1r^{-1}dV(r))]= \exp (t (\log \phi_{V(1)}(t))^\prime), \ \ t \in\Rset.
\end{equation}

\medskip
As an application of  Proposition 1 we find three expressions for the BDDF of  the gamma $\gamma_{\alpha,\lambda}$ distribution, i.e., the probability distribution with the  density $\lambda^{\alpha}/\Gamma(\alpha)\,x^{\alpha-1}e^{-\lambda x}1_{(0,\infty)}(x)$, ($\alpha>0, \lambda>0)$. Recall that,
\begin{multline}
\phi_{\gamma_{\alpha, \lambda}}(t)=\mathbb{E}[e^{it\gamma_{\alpha, \lambda}}]=(1-it/\lambda)^{-\alpha}=\exp[\alpha \int_0^\infty (e^{itx}-1)\,\frac{e^{-\lambda x}}{x}dx]
\\ = \exp\int_0^1 \alpha\, [\int_0^\infty (e^{itxv} -1) \, \lambda e^{-\lambda x} dx]\,\frac{dv}{v},  \  (\phi_{\gamma_{\alpha,\lambda}} \ \mbox{in terms of its BDCF)}
\end{multline}
where the last equality follows from the random integral representation of the gamma distribution;  cf. Jurek (1996), Corollary 1 and Remark 1.
\begin{lem}
For a gamma variable $\gamma_{\alpha,\lambda}$ its BDRV  equals  $Y(1)=  \sum_{k=1}^{N_{\alpha}}\mathcal{E}_k (\lambda)$, where $N_{\alpha}, \ \mathcal{E}_k(\lambda), k=1,2,...$  are independent, $N_{\alpha}$  is Poisson variable with a parameter $\alpha$  and $\mathcal{E}_k(\lambda)$ are exponential  identically distributed with parameter $\lambda$.  Moreover,  the  BDDF  for the selfdecomposable $X= \gamma_{\alpha,\lambda}$  is  given as follows:
\begin{multline}
G_{X}(a)=P\big( \sum_{k=1}^{N_{\alpha}}\mathcal{E}_k (\lambda) \le a \big)\\ = \frac{1}{2}+\frac{1}{\pi}\int_0^\infty \exp(-\frac{\alpha\,(t/\lambda)^2}{1+ (t/\lambda)^2})\, \sin \big(t a- \frac{\alpha\,(t/\lambda)}{1+(t/\lambda)^2}\big)\,\frac{1}{t}dt \\
 = e^{-\alpha}+ e^{-\alpha} \int_0^{2\sqrt{\alpha  \lambda a}}I_1(w)e^{-w^2/4\alpha}dw, \  \ \ a\in C_{G_X}\cap (0,\infty);  \qquad \qquad
\end{multline}
\mbox{ where $I_1(x)$ is a Bessel function; cf. Appendix c) below.}
\end{lem}
 More examples illustrating Proposition 1 are  given in Section III, below.

\medskip
\textbf{I. 2.} \emph{ Series representations for log-gamma variables.}

\noindent From Shanbhag, Pestana and Sreehari (1977), Remark 1, p.291,  we know that  $\log\gamma_{\alpha, \lambda}$  random variables  are  selfdecomposable, so  they have the representation (1).  However, they  also  admit the following series representation:
\begin{prop}
Let $\mathbf{C}$  be the Euler constant and $ \gamma_{1, \alpha+n}, n=0,1,2,...$ be a sequence of independent exponential random variables with scale parameters $\alpha +n$. Then the series
\begin{multline}
S: =  -\log \lambda - \mathbf{C} - \gamma_{1, \alpha}- \sum_{n=1}^\infty(\gamma_{1, \alpha +n} - 1/n) \\
= -\log \lambda +\Psi(\alpha)-\sum_{n=0}^\infty(\gamma_{1, \alpha +n} - 1/(\alpha+n)),  \ \ \ \
\end{multline}
converges almost surely (in probability, in distribution, in $ L^2$).
Moreover, $S\stackrel{d}{=}\log \gamma_{\alpha, \lambda}$ and  $ \phi_{\log\gamma_{\alpha,\lambda}}(t)=
\lambda^{- it}\frac{\Gamma(\alpha+it)}{\Gamma(\alpha)}$. \ ($\Gamma$ \mbox{denotes the Gamma function}).
\end{prop}
From the above series representation we conclude
\begin{cor}
For log- gamma variables we have
\[
\mathbb{E}[\log\gamma_{\alpha,\lambda}]= - \log\lambda + \Psi^{(0)}(\alpha); \ \  Var[\log\gamma_{\alpha,\lambda}]=\Psi^{(1)}(\alpha)=\int_0^\infty \frac{x e^{- \alpha x}}{1-e^{-x}}dx.
\]
($\Psi^{(0)}$ and $\Psi^{(1)}$ denote the digamma and its first derivative function, respectively; cf. Gradsteyn and Ryzhik (1994), Sec.8.36.
\end{cor}
Our next aim is to find the BDCF $\psi_{\log\gamma_{\alpha, \lambda}}(t)$ for  the selfdecomposable characteristic function $\phi_{\log\gamma_{\alpha, \lambda}}(t)$.
\begin{prop}
The background driving characteristic function (BDCF) for the log $\gamma_{\alpha,\lambda}$ variable is as follows
\begin{equation*}
\psi_{\log\gamma_{\alpha, \lambda}}(t) = \exp[\,it\,(-\log\lambda+\Psi^{(0)}(\alpha))+\int_0^\infty(e^{-itx} - 1  + itx) e^{-\alpha x}\,h_{\alpha}(x) dx],
\end{equation*}
where $h_{\alpha}(x):= [\alpha +(1-\alpha)e^{-x}](1-e^{-x})^{-2}$  and  $\int_0^\infty x^2 e^{-\alpha x}h_{\alpha}(x)dx<\infty$ and $\Psi^{((0)}(z)$ is the digamma function.
\end{prop}
Since $\log\gamma_{\alpha, \lambda}$  is selfdecomposable, it has a background driving distribution function  (BDDF) $G_{\log \gamma_{\alpha, \lambda}}$ for which we have:
\begin{cor}
 Let $G_{\log\gamma_{\alpha, \lambda}}$ be the BDDF of the log-gamma varaible.  Then its expected value and the variance are  as follows:

\noindent (a) \ $\mathbb{E}[G_{\log\gamma_{\alpha, \lambda}}(x)]= -\log\lambda +\Psi^{(0)}(\alpha)$;

\noindent (b) \  $Var[ G_{\log\gamma_{\alpha,\lambda}}]=\int_0^\infty x^2e^{-\alpha x}[\alpha +(1-\alpha)e^{-x}](1-e^{-x})^{-2}dx$.

\noindent (c) \  For $\beta>0$ \, all moments  $\mathbb{E}[\, |x|^{\beta}dG_{\log\gamma_{\alpha, \lambda}}(x)]$ are finite.
\end{cor}
Finally, here is a series  representation for the innovation variable  $X_c$  from (2)  for the log-gamma variable:
\begin{prop}
For independent binomial
$b_c^{(k)}$ and gamma  random variables $\gamma_{1,\alpha +k}, \  k=0,1,2,...$ where $P(b_c^{(k)}=1)=1- c $ and $P(b_c^{(k)}=0)=c $  with  $0<c<1$ ,  the random series
\begin{equation}
Z_c(\alpha):= - (1-c) \mathbf{C} + b_c^{(0)}(-1) \gamma_{1,\alpha} +\sum_{n=1}^\infty[ b_c^{(n)}(-1)\gamma_{1,\alpha +n} +  \frac{1-c}{n}],
\end{equation}
converges almost surely (in probability, in distribution). Moreover,  $Z_c(\alpha)$  has characteristic function
$\mathbb{E}[e^{itZ_c(\alpha)}]= \lambda^{-it(1-c)}\frac{\Gamma(\alpha+it)}{\Gamma(\alpha+ict)} = \phi_{\log\gamma_{\alpha,\lambda}}(t)/  \phi_{c \log\gamma_{\alpha,\lambda}}(t)$.
\end{prop}
\textbf{II. PROOFS.}

\emph{Proof of Proposition 1.}

\noindent If $\phi$ and $\psi$ are characteristic functions of $X$ and $Y(1)$, respectively, then
\begin{equation}
 \phi(t)= \exp\int_0^t \log\psi(u)\frac{du}{u},  \ \mbox{equivalently} \ \ \psi(t)=\exp [\,t\, \frac{\phi^\prime (t)}{\phi(t)}\,] \ \  t\neq 0;
\end{equation}
cf. Jurek (2001), Proposition 3.  If  $Y(1)$ is the BDRV of selfdecomposable random variable X then $Y(1)$ has finite logarithmic moment; cf. Jurek and Vervaat (1983), Theorem 2.3.  Consequently, by  Gil-Pelaez (1951), Wendel (1961) and Ushakov (1999), Theorem 1.2.4 (cf. also Appendix $(b))$) we have that the cumulative distribution function of $Y(1)$ can be obtained from $\psi$ via the following inversion formula:
\[
P(Y(1)\le a)=\frac{1}{2}-\frac{1}{\pi}\int_0^\infty \Im( e^{-ita} \cdot \psi(t))\,\frac{dt}{t},  \ \   (\mbox{$a$ is a continuity point});
\]
 which with (10) completes the argument for the proof of Proposition 1.

\medskip
\emph{Proof of Lemma 1.}

From (6), $\phi_{\gamma(\alpha,\lambda)}(t)= (1-it/\lambda)^{-\alpha}$ and therefore from (10) we find that their  BDCF are  given as follows
\begin{equation}
\psi_{\gamma(\alpha,\lambda)}(t)=\exp[\alpha\, \frac{ i t/\lambda}{1-it/\lambda}]=  \exp\alpha\,[ \frac{1}{1-it/\lambda}-1].
\end{equation}
Thus they correspond to compound Poisson distributions. More explicitly , if $\mathcal{E}_k(\lambda), k=1,2,..$ are i.i.d  (exponentially distributed with scale parameter $\lambda$) and independent of $N_{\alpha}$, Poisson variable with parameter $\alpha$,  then  for
\begin{equation}
Y(1) := \sum_{k=1}^{N_{\alpha}}\mathcal{E}_k (\lambda)\ \ \ \mbox{we have} \ \ \mathbf{E}[e^{i t Y(1)}]= \exp\alpha\,[ \frac{1}{1-it/\lambda}-1],
\end{equation}
that is,  $Y(1)$, has a compound Poisson distribution and it  is the BDRV for the selfdecomposable $\gamma_{\alpha,\lambda}$ random variable. ( If $N_\alpha=0$  then  the  sum  in (12) is  zero.) This completes a proof of the  first equality in (7).

For the second equality in (7) we use Proposition 1. Namely, for a>0,
\begin{multline}
P\big( \sum_{k=1}^{N_{\alpha}}\mathcal{E}_k (\lambda) \le a \big)=\frac{1}{2}-\frac{1}{\pi}\int_0^\infty \Im(\exp(-ita +  \alpha \frac{it/\lambda}{1- it/\lambda}))\,\frac{dt}{t} \\ =
 \frac{1}{2}+\frac{1}{\pi}\int_0^\infty \exp(-\frac{\alpha\,(t/\lambda)^2}{1+ (t/\lambda)^2})\, \sin \big(t a- \frac{\alpha\,(t/\lambda)}{1+(t/\lambda)^2}\big)\,\frac{1}{t}dt.
\end{multline}

For the third equality in (7) we use  the additive property of gamma function with respect to the shape parameter  and the identity
\begin{equation}
\sum_{k=1}^\infty \frac{b^{k-1}}{k! (k-1)!}=\frac{I_1(2\sqrt{b})}{\sqrt{b}},
\end{equation}
(cf. Gradshteyn-Rhyzik(1994), \textbf{8.445} or Appendix below),
as follows
\begin{multline*}
P\big( \sum_{k=1}^{N_{\alpha}}\mathcal{E}_k (\lambda) \le a \big)=
\sum_{k=0}^\infty P[( \sum_{j=1}^{N_{\alpha}}\mathcal{E}_j (\lambda) \le a)
\cap (N_\alpha = k) \big ] \\ = e^{-\alpha} + \sum_{k=1}^\infty  e^{-\alpha}
\frac{\alpha^k}{k!} P(\gamma_{k,\lambda}\le a)=
e^{-\alpha} + e^{-\alpha}  \sum_{k=1}^\infty \frac{\alpha^k}{k!} \frac{\lambda^k}{\Gamma(k)}\int_0^a x^{k-1}e^{-\lambda x}dx \\
e^{-\alpha}+e^{-\alpha} \alpha\lambda \int_0^\alpha \sum_{k=1}^\infty \frac{(\alpha \lambda x)^{k-1}}{k! (k-1)!}dx =  e^{-\alpha}+ e^{-\alpha}\alpha\lambda \int_0^a\frac{I_1(2\sqrt{\alpha \lambda x})}{\sqrt{\alpha\lambda x}} e^{-\lambda x} dx,
\end{multline*}
and changing variable we conclude the proof of  Lemma 1.

\medskip
\emph{Proof of Proposition 2.}

Firstly, recall that  series of centered  independent gamma variables

\noindent $\sum_{n=1}^\infty ( \gamma_{1,\alpha+n} - 1/(\alpha+n))$ converges in all above mentioned modes, because $ \inf_n(n+\alpha)=\alpha +1>0$ and $\sum_n(n+\alpha)^{-2}<\infty$ ; cf. Jurek (2000), Proposition 1 and Corollary 2.

Secondly, we have  two converging  numerical series
\begin{equation}
\sum_{n=1}^\infty \frac{1}{n(n+\alpha)}= \alpha^{-1}(\Psi^{(0)}(\alpha+1)+\mathbf{C}); \  \  \ \sum_{n=1}^\infty \frac{1}{(n+\alpha)^2}=\Psi^{(1)}(\alpha+1);
\end{equation}
 cf. Appendix part \textbf{c).}   Above $\Psi(z)\equiv \Psi^{(0)}(z):= \frac{d}{dz}\log \Gamma(z)$ is the digamma function and $\Psi^{(1)}(z)$ is  its first derivative .

 Finally,  since $ \gamma_{1,\alpha+n}- 1/n =   \gamma_{1,\alpha+n} - \mathbb{E}[\gamma_{1,\alpha+n}] -\alpha/n(\alpha+n)$  we conclude that
the infinite random series in (8) converges in all three  modes.

Since, in particular, series (8)  converges  weakly  its characteristic function  is given as an infinite product :
\begin{multline}
 \mathbb{E}[e^{it S}]=\lambda^{-it}\,\, \exp(-it \mathbf{C}) \,\, (1+it/\alpha)^{-1}\,\,  \prod_{n=1}^\infty e^{it/n}(1+it/(\alpha+n))^{-1}  \ \ \  \mbox{(by (6))} \\
= \lambda^{-it} \exp\Big[-it\mathbf{C} + \int_0^\infty (e^{-itx}-1)\frac{e^{-\alpha x}}{x}dx  \\ + \sum_{n=1}^\infty \Big( it/n+\int_0^\infty(e^{-itx}-1)\frac{e^{-x(\alpha+n)}}{x}dx)\Big) \Big] \\ =  \lambda^{-it} \exp\Big[-it\mathbf{C} - it \alpha^{-1} + \int_0^\infty (e^{-itx}-1 + i tx)\frac{e^{- \alpha x}}{x}dx \\ +\sum_{n=1}^\infty \Big( it (\frac{1}{n} - \frac{1} {\alpha+n}
)+\int_0^\infty(e^{-itx}-1  + it x)\frac{e^{-x(\alpha+n)}}{x}dx)\Big) \Big]  \  \  (\mbox{by (15)}) \\ =  \lambda^{-it} \exp\Big[-it\mathbf{C}   + it(\Psi^{(0)}(\alpha+1) - \alpha^{-1} +\mathbf{C}) \\ + \int_0^\infty(e^{-itx}-1  + it x)(\sum_{n=0}^\infty\frac{e^{-x(\alpha+n)}}{x})dx)\Big]  \\
=    \lambda^{-it} \exp[ it \Psi^{(0)}(\alpha) +    \int_0^\infty(e^{-itx}-1  + it x) \frac{e^{-\alpha x}}{x(1-e^{-x})}dx  ] \\ =\lambda^{-it}\,\frac{\Gamma(\alpha +it)}{\Gamma(\alpha)} = \mathbb{E}[e^{it \log\gamma_{\alpha, \lambda}} ]. \qquad \qquad \qquad \qquad \quad \qquad \qquad
\end{multline}
The last two equalities are from Jurek (1997), example (c) on  p. 98, where the correct reference  should be to Whittaker and Watson (1920) p. 249 and/or Shanbhag, Pestana and Sreehari (1977)), Lemma 1. This gives a proof of the first equality in (8).

For the second equality in (8), note that by (15) and Appendix \textbf{(c)},(iv) we have
$$
\sum_{n=1}^\infty(\gamma_{1,\alpha+n}-\frac{1}{n})=\sum_{n=1}^\infty(\gamma_{1,\alpha+n}-\frac{1}{n+\alpha})-\Psi(\alpha)- 1/\alpha-C,
$$
 which concludes the proof of Proposition 2.

[In Bondesson (1992), Chapter 7, p. 113 there is a series representation of $\log\gamma_{\alpha, 1}$ as in the second part of (8), without specifying the constant. More importantly, it was obtained as an example of analytically introduced class of extended generalized Gamma convolution (EGGC); cf. Definition on p. 103.]
\begin{rem}  (a) \emph{
 For $\alpha>0$ and $t\in\Rset$ we have
\begin{multline}
 \frac{ e^{-i \mathbf{C}t}}{1+i t/\alpha}\prod_{n=1}^\infty \frac{\exp(i t/n)}{1+\frac{it}{\alpha+n}} =
\frac{\Gamma(\alpha+it)}{\Gamma(\alpha)}=  \\ =  \exp[ it \Psi^{(0)}(\alpha) +    \int_{- \infty}^0(e^{itx}-1  - it x) \frac{e^{\alpha x}}{|x|(1-e^{x})}dx ].
\end{multline}
(b) \  \   $\mathbb{E}[(\exp it\log\gamma_{\alpha, 1})]=\frac{\Gamma(\alpha+it)}{\Gamma(\alpha)}, \ t\in\Rset.$}
\end{rem}
\noindent For part (a) : the first equality coincides with GR \textbf{8.326}(2) and with Erdelyi (1953), p. 5, formula (3).  The second one is from  Shanbhag, Pestana and Sreehari (1977), Lemma 1; note that their constant $c_{\alpha}=\Psi^{(0)}(\alpha)$ via (16).

\noindent  Part (b) is from last line in (16).

\emph{Proof of Corollary 2.}

Using the representation (8) and formula (15) we have
\[
\mathbb{E}[\log\gamma_{\alpha,\lambda}]=- \log \lambda - \mathbf{C}-\alpha^{-1}+ \alpha\sum_{n=1}^\infty(\frac{1}{n(\alpha +n)} = - \log \lambda +\Psi^{(0)}(\alpha +1) - \alpha^{-1}
\]
\[= - \log\lambda + \Psi^{(0)}(\alpha); \  [\mbox{here we used identity} \  \Psi^{(0)}(\alpha +1)-\Psi^{(0)}(\alpha)= 1/\alpha]
\]
which gives the first identity.  For the second one,  using independence of the summands in (8) and the equality
\noindent $\sum_{n=1}^\infty\frac{1}{(\alpha+n)^2}= \Psi^{(1)}(\alpha +1)$
we  get
$
Var[\log\gamma_{\alpha,\lambda}]=1/\alpha^2 + \sum_{n=1}^\infty\frac{1}{(\alpha +n)^2}= \Psi^{(1)}(\alpha +1) +1/\alpha^2=\Psi^{(1)}(\alpha),$  which completes the first formula for the variance.

\noindent Second one follows from (16) viewed as Kolmogorov  representation of infinitely divisible variables with finite second moments. This completes the proof of Corollary 2.

\emph{Proof of Proposition 3.}

From (8) and first line in (6)  we have
\begin{multline*}
\phi_{\log \gamma_{\alpha, \lambda}}(t)= e^{-it(\mathbf{C}+\log \lambda)}\,\mathbb{E}[e^{-it \gamma_{1,\alpha}}] \,
\prod_{n=1}^\infty\mathbb{E}[e^{-it(\gamma_{1,\alpha+n} -1/n)}]\\ =
 e^{-it(\mathbf{C}+\log \lambda)} (1+it/\alpha)^{-1}\ \prod_{n=1}^\infty e^{it/n} (1+ it/(\alpha +n))^{-1}.
\end{multline*}
From second line in (6), for $\beta>0$, we have identity
\[
(1+it/\beta)^{-1}=\exp\Big[ -it/\beta + \int_0^1  [\int_0^\infty(e^{-itux}-1 +itux)\beta e^{-\beta x})dx]\frac{du}{u}   \Big]
\]
Applying that identity for $\beta:= \alpha+n, n=0,1...$ in the first formula above and then using (15)  we
get
\begin{multline}
\lambda^{it} \phi_{\log\gamma_{\alpha, \lambda}}(t)=  \exp\Big[-it\mathbf{C} - it \alpha^{-1} + \int_0^1[\int_0^\infty (e^{-itux}-1 + i tux)\alpha e^{- \alpha x}dx]\frac{du}{u} \\ +\sum_{n=1}^\infty \Big( it (\frac{1}{n} - \frac{1} {\alpha+n}
)+\int_0^1[\int_0^\infty(e^{-itux}-1  + itu x) (\alpha+n)e^{-x(\alpha+n)}dx] \frac{du}{u}\Big)\Big]  \\
=\exp\Big[ it \Psi^{(0)}(\alpha) + \int_0^1[\int_0^\infty(e^{-itux}- 1  + itu x)\sum_{n=0}^\infty (\alpha+n)e^{-x(\alpha+n)}dx] \frac{du}{u}\Big ]\\
= \exp \int_0^1 [ itu \Phi^{(0)}(\alpha)+\int_0^\infty(e^{-itux} - 1  + itu x) e^{-\alpha x}(\frac{\alpha}{1-e^{-x}} +
\frac{e^{-x}}{(1-e^{-x})^2})dx]\frac{du}{u}\\
=\exp \int_0^1 \big[ itu \Phi^{(0)}(\alpha)+\int_0^\infty(e^{-itux} - 1  + itu x) e^{-\alpha x}\,h_{\alpha}(x) dx\big]\frac{du}{u}.
\end{multline}
From (10), the relations between selfdecomposable $\phi$ and its BDCF $\psi$, are of the form
$\phi(t)=\exp\int_0^1\log\psi(tw)\frac{dw}{w}$. Hence we infer the BDCF $\psi_{\log\gamma_{\alpha, \lambda}}(t) $ is in the square bracket in (18)  by putting u=1.

Since for $\alpha>0$ and $x>0$
\begin{multline*}
\alpha +(1-\alpha)e^{-x}\le \max (\alpha, 1-\alpha) + \max(\alpha,1-\alpha) e^{-x}\le  2 \max(\alpha, 1-\alpha);\\
(\frac{x}{1-e^{-x}})^2\le (1+x)^2 \ \mbox{ and} \ \int_0^\infty e^{-\alpha x}(1+x)^2dx<\infty \ \mbox{ because} \\  \int (1+x)^2e^{-\alpha x}dx= -\alpha^{-3}\,e^{-\alpha x}[(\alpha(x+1)+1)^2 +1] +const;
\end{multline*}
 we infer that the integral is finite, which completes the proof of Proposition 3.

\medskip
\emph{Proof of Corollary 3.}

\noindent Parts (a) and (b) are from the second to last line in  (16) and the  Kolmogorov's representation of infinitely divisible variables with finite second moments. Part (c) follows from the fact that $\int_1^\infty x^{\beta}e^{-\alpha x}h_{\alpha}dx <\infty$ and from Jurek and Mason(1993), Proposition 1.8.13 on p. 36.

\medskip
\emph{Proof of Proposition 4.}

\noindent\emph{Step 1.}
Let us put $\xi_n:=[ b_c^{(n)}(-1)\gamma_{1,\alpha +n} +  \frac{1-c}{n}]$. Then
\[
\mathbb{E}[\xi_n]=  - \frac{1-c}{\alpha +n}+ \frac{1-c}{n}=\frac{\alpha(1-c)}{n(\alpha+n)}, \ \ \sum_{n=1}^\infty \mathbb{E}[\xi_n]= (1-c)(\Psi^{(0)}(\alpha+1)+\mathbf{C});
\]
(see Appendix, Section \textbf{c).} or use WolframAlpha)
so, series of expected values of $\xi_n$ converges. Moreover, note that
\begin{multline*}
Var[\xi_n]=\mathbb{E}[( b_c^{(n)}(-1)\gamma_{1,\alpha +n})^2]- \big(\mathbb{E}[ b_c^{(n)}(-1)\gamma_{1,\alpha +n}]\big)^2=  \frac{1-c}{(\alpha+n)^2}-    \frac{(1-c)^2}{(\alpha+n)^2}\\= \frac{c(1-c)}{(\alpha+n)};   \qquad   \sum_{n=1}^\infty Var[\xi_n]= c(1-c)\sum_{n=1}^\infty \frac{1}{(\alpha+n)^2}= c(1-c)\Psi(\alpha+1)
\end{multline*}
so the series of  variances converges as well. All in all,  Three Series Theorem gives convergence almost surely
(and in probability) of the series in (9).

\noindent\emph{Step 2.}  Since for  independent variables  binomial $b_c$ ($P(b_c=1)= 1-c)$ and  $\gamma_{1,\beta}$ we have that
\[
\mathbb{E}[e^{it b_c\gamma_{1, \beta}}]= \mathbb{E}[\mathbb{E}[e^{itb_c\gamma_{1,\beta}}| b_c]= c + \frac{1-c}{1-it/\beta}=\frac{1-i ct/\beta}{1- it/\beta},
\]
therefore for $\xi_n$ defined above
\begin{multline}
\mathbb{E}[e^{it \xi_n}]= e^{it(1-c)/n} [ c + \frac{1-c}{1+ it/(\alpha+n)}]= e^{it(1-c)/n} \frac{1 +i ct/(\alpha +n)}{1+it/(\alpha+n)} \\ =
\big[ e^{it/n}\,\frac{1}{1+it/(\alpha +n)} \big]\cdot \big[ e^{ict/n}\,\frac{1}{1+ict/(\alpha +n)} \big]^{-1}.
\end{multline}

\noindent\emph{Step 3.} Using the definition  of $Z_c(\alpha)$ and then (19) and Remark 2  we have
\begin{multline*}
\mathbb{E}[e^{it Z_c(\alpha)}]=    e^{-it \mathbf{C}(1-c)}\mathbb{E}[e^{itb_c(-1)\gamma_{1,\alpha}}]\prod_{n=1}^\infty\mathbb{E}[e^{it \xi_n}]  \\
= e^{-it \mathbf{C}(1-c)}\, \frac{1 +itc/\alpha}{1+it/\alpha} \ \prod_{n=1}^\infty e^{it(1-c)/n} \  \frac{1 +itc/(\alpha +n)}{1+it/(\alpha+n)}  \\ =  [ \Gamma(\alpha+it)/\Gamma(\alpha)]/[\Gamma( \alpha +i ct)/\Gamma(\alpha)]=  \mathbb{E}[e^{it \log\gamma_{\alpha,\lambda}}]/\mathbb{E}[e^{it c\log\gamma_{\alpha,\lambda}}];
\end{multline*}
which completes the proof.

\medskip
\textbf{III. MORE EXAMPLES AND NEW FORMULAS.}

\textbf{1.} \ \emph{L\'evy distributions $L(m,c)$.}

For  a location parameter $ m\in \Rset$ and a scale parameter $c>0$, L\'evy variable L\'evy(m,c) has the probability density function
\begin{equation}
f(m, c; x):=\sqrt{\frac{c}{2\pi}}\,\,\frac{ e^{-\frac{c}{2(x-m)}}}{(x-m)^{3/2}} , \ \  \ \mbox{for} \ \ x > m.
\end{equation}
If one defines the error function $erf(x):= \frac{2}{\sqrt{\pi}}\int_0^x e^{-s^2}ds$, ($erf(\infty)=1$) then the cumulative probability  distribution function $F\equiv F(m,c;x)$ of $L(m,c)$  is given as
\begin{multline}
F(m, c;  x)= erfc \big( \sqrt{\frac{c}{2(x-m)}}\big ):= 1-erf \big( \sqrt{\frac{c}{2(x-m)}}\big ) \\= \frac{2}{\sqrt{\pi}}\int_{ \sqrt{\frac{c}{2(x-m)}}}^\infty e^{-s^2}ds, \   \qquad  \mbox{for} \ x>m. \qquad \qquad \qquad
\end{multline}
Finally the  characteristic function $\phi_F$ is equal to
\begin{equation}
\phi_F(t)= e^{i m t-|ct|^{1/2}(1- i\,sign\, t) },  \ \ \ t\in\Rset.
\end{equation}
Hence F= L\'evy(m,c) are the  stable distributions with the exponent 1/2. Consequently, they are selfdecomposable and thus admit the  random integral representation (1). From (22), via (10), we get
\begin{equation}
\psi_F(t)= \exp[i m t-|(c/2)t|^{1/2}(1- i\,sign\,t)], \mbox{i.e.},   Y(1) \stackrel{d}{=} \mbox{ L\'evy (m, c/2)}.
\end{equation}
So, the BDPD for 1/2-stable $L(m,c)$  is another 1/2-stable  distribution but  with the  scale parameter $c/2$, that is, $L(m,c/2)$.

\medskip
From Proposition 1, the cumulative probability distribution function G of Y(1) is
\begin{multline}
G(x)= \frac{1}{2}- \frac{1}{\pi}\int_0^\infty\Im[\exp(- itx + i m t - ( (c/2) t)^{1/2} (1- i)]\frac{dt}{t} \\
= \frac{1}{2}+\frac{1}{\pi}\int_0^\infty\ e^{- ((c/2)t)^{1/2}}\sin (t (x-m) - ((c/2)t)^{1/2})\frac{dt}{t} \ \ ( s:= ((c/2)t)^{1/2})\\
= \frac{1}{2}+\frac{2}{\pi}\int_0^\infty\ e^{- s}\sin( (2/c)(x-m) s^2- s)\frac{ds}{s} = erfc \big(  \sqrt{\frac{c/2}{2(x-m)}}  \big),
\end{multline}
where the last equality is justified by (21).

As a by-product of the above we get the identity:
\begin{cor}
For $a\in\Rset$ we have identity
$$
\int_0^\infty e^{-x}\sin((a^2x^2-x)\frac{dx}{x}= \frac{\pi}{2} \Big( erfc ( (|a| \sqrt{2})^{-1})- \frac{1}{2}\Big)=
 \frac{\pi}{2} \Big( \frac{1}{2}- erf( (|a| \sqrt{2})^{-1}))
$$
\end{cor}
\textbf{2.} \emph{Stable distributions.}
The above straightforward calculations showed that the BDDF of 1/2-stable distribution is 1/2-stable distribution.
That fact is  not surprising one as the stable laws are the  fixed points of the integral mapping (1); cf. Jurek and Vervaat (1983), Theorem 4.1.

\textbf{3.} \emph{Bessel-K distributions}

\noindent Bessel-K   distribution is the symmetrization of  $\gamma_{\alpha, \lambda}$, so it is selfdecomposable and its characteristic function is
$\phi_{BK(\alpha, \lambda)}(t)= (1+t^2/\lambda^2)^{-\alpha}$. \  [$ K$  stands here for the Bessel function  $K_{\alpha}$ which appears in the formula for the probability density function of Bessel-K distributions; cf. Johnson, Kotz and Balakrishnan(1994), Sect. 4.4,  p. 50].   Moreover, we have
\begin{equation}
BK(\alpha, \lambda)=\sqrt{2 \gamma_{\alpha,\lambda^2}}\, \, Z,
\end{equation}
where the gamma variable and the standard normal  $Z$ are independent. This is so, because
 \begin{multline*}
\mathbb{E}[e^{it \sqrt{2 \gamma_{\alpha, \lambda^2}}Z}]=\mathbb{E}\big[\mathbb{E}[e^{it \sqrt{2 \gamma_{\alpha, \lambda^2}}Z }| \gamma_{\alpha,\lambda^2}]\big] \\ =\mathbb{E}[e^{-t^2\gamma_{\alpha,\lambda^2}}]   = \frac{\lambda^{2 \alpha}}{\Gamma(\alpha)}\int_0^\infty e^{-t^2 x}x^{\alpha-1}e^{-\lambda^2 x}dx  =  \frac{\lambda^{2 \alpha}}{\Gamma(\alpha)}\int_0^\infty e^{- (t^2+ \lambda^2)x}x^{\alpha-1}dx \\= [\frac{\lambda^2}{\lambda^2+t^2}]^\alpha  = [\frac{1}{1+(t/\lambda)^2}]^\alpha = \phi_{BK(\alpha,\lambda)}(t)  \qquad \qquad \qquad
\end{multline*}
which concludes the proof of (26).

Hence and from (10)  the  BDCF  of  $BK(\alpha, \lambda)$ is equal to
\begin{equation}
\mathbb{E}[e^{itY(1)}]=\psi_{BK(\alpha,\lambda)}(t)=\exp[ 2\alpha(\frac{1}{1+t^2/\lambda^2}-1)].
\end{equation}
 \noindent Consequently, the  BDLP $(Y(t),t\ge 0)$ for Bessel-K  distribution  is   the compound Poisson process.
\begin{cor}
If a random variable $BK$ has the Bessel-K distribution then its innovation $BK_c$ has the following representations:
\[
\mathbb{E}[\exp(itBK_c)] =[c^2 1+(1-c^2)\frac{1}{1+t^2/\lambda^2}]^\alpha \stackrel{d}{=}(b_{c^2}\mathcal{E}
(\lambda^{\circ}))^{\star\alpha}\stackrel{d}{=}\sum_{k=1}^{N_{- \alpha \log c^2}}c^{\eta_k}\mathcal{E}(\lambda^{\circ}_k),
\]
where $\eta_k$  are i.i.d.  uniform on $(0,1)$,  $\mathcal{E}(\lambda^{\circ}_k)$  are i.i.d. symmetric  exponential variables
and $\stackrel{d}{=}$ means the equality in distribution.
\end{cor}
\noindent \emph{Proof.}  Since  $BK_c:= \int_0^{-\log c } e^{-s}dY(s)$, with $Y(1)$ given in (27), therefore it can  be written as  follows
\begin{multline*}
\mathbb{E}[\exp(it\,BK_c)]= \exp\int_0^c \log\psi_{BK(\alpha,\lambda)}(e^{-s}t)ds \\=
\exp\int_0^{-\log c}( -2 \alpha)\,\frac{t^2 e^{-2s}/\lambda^2}{1+t^2e^{-2s}/\lambda^2}ds= (\frac{1+c^2t^2/\lambda^2}{1+t^2/\lambda^2})^\alpha \\ =[c^2 1+(1-c^2)\frac{1}{1+t^2/\lambda^2}]^\alpha \stackrel{d}{=}(b_{c^2}\mathcal{E}(\lambda^{\circ}))^{\ast \alpha}\stackrel{d}{=}\sum_{k=1}^{N_{- \alpha \log c^2}}c^{\eta_k}\mathcal{E}(\lambda^{\circ}_k),
\end{multline*}
To see the last two equality  one needs to compute the characteristic functions of $b_c^2\,\mathcal{E}(\lambda)$ and of the random sum (compound Poisson distribution).

\textbf{Note:} analogous formulas are valid for $\gamma_{\alpha, \lambda}$ variables. Also
compare  Lawrance (1982) and Hosseini (2019).

\medskip
\textbf{IV. APPENDIX.}

\noindent $\textbf{a).}$  Besides those two equivalent characterizations  (formulas (1) and (2)) of the selfdecomposability, also  the  following  limit laws  give the class  L variables $X$.

Namely for $X\in L$, there exist \underline{a strong mixing sequence}  $(V_n, n=1,2,...) $ of random variables, there exist deterministic sequences $a_n>0, b_n \in \Rset, n=1,2,...$ such that

(i) the triangular array $(a_nV_j: 1\le j \le n; n\ge 1)$ is infinitesimal; and

(ii) $a_n(V_1+V_2+...+V_n) +b_n \Rightarrow X,  \ n\to \infty$ ( weak convergence);

\noindent cf. Bradley and Jurek (2016). [Of course, sequences of independent variables are strongly mixing.]

\medskip
\noindent $\textbf{b).}$ \  For  the proof of  Proposition 1 the following inversion formula was essential.

\emph{ Let distribution function G, with its characteristic function $\psi$, has finite logarithmic moment, that is, $\int_{-\infty}^\infty \log(1+|x|)dG(x)< \infty$. Then}
\begin{equation}
G(x)=\frac{1}{2}-\frac{1}{\pi}\int_0^\infty \Im(e^{-iux}\psi(u))\frac{du}{u}, \ \ \ x \in C_G
\end{equation}
cf. Gil-Paleaz(1951), Wendel(1961) and Ushakov(1999), Theorem 1.2.4.

\medskip
\noindent $\textbf{c).}$  For an ease of reference let us  recall some special functions and formulas:
\begin{multline*}
 (i) \  \Gamma(z):= \int_0^\infty x^{z-1}e^{-x}dx, \   \Re z >0; \
(ii) \ \ \Psi(z):=\frac{d}{dz}\log\Gamma(z) ,  \ \Re z>0;
\\
(iii) \ \   I_{\alpha}(z)= \sum_{j=0}^\infty ,\frac{ (z/2)^{\alpha + 2j}}{j!\Gamma(\alpha+j+1)}, \  \alpha>0.\ (\mbox{cf.\textbf{8.445}}); \ \ \qquad \qquad
\\
(iv) \ \  \Psi(x)=-\textbf{C}-\frac{1}{x}+x\sum_{k=1}^\infty\frac{1}{k(k+x)}; \ \   (\mbox{cf.\ \textbf{8.362}.1}); \\  (v)\ \ \Psi(x+1)= \Psi(x)+ \frac{1}{x}; \ \ \ (\mbox{cf.\ \textbf{8.365}.1});
\\
(vi)  \ \Psi^{(n)}(x)= (-1)^n n!\sum_{k=0}^\infty\frac{1}{(x+k)^{n+1}}, \ \  n\ge 1; \ \  (\mbox{cf.\ \textbf{8.363}.8;}) \\ \mbox{(All bold face formulas are taken from Gradshteyn-Rhyzik(1994)).}
\end{multline*}
\textbf{Acknowledgments.}
(1) The author would like to thank Andrew Stuart and Bamdad Hosseini for their kind invitation to the California Institute of Technology, and for an introduction to the problems of algorithms for simulating random variables (in particular, for selfdecomposable ones). Some of the above results were inspired by conversations at Caltech.

(2) Referee's comments and suggestions helped to improve the language and presentation of the results.

\medskip
REFERENCES

[1] L. Bondesson (1992), \emph{Generalized gamma convolutions and related classes of distributions and densities,} Lectures Notes in Statistics, vol.76; Springer-Veralg, New York

[2] R. C. Bradley. Z. J. Jurek (2016),  Strong mixing and  operator- selfdecomposability, \emph{J. Theor. Probability}, \textbf{29},  2016, pp. 292-306.

[3] A. Erdelyi, W. Magnus, F.Oberhettinger, F.G. Tricomi (1953), \emph{Higher transcendental functions}, vol. I,  New York, Toronto, London: McGraw-Hill Book Company.

[4] W. Feller (1966), \emph{An introduction to probability theory and its applications}, vol II, John Wiley $\&$ Sons,Inc. New York, London, Sydney.

[5] J. Gil-Pelaez (1951),  Note on the inversion theorem, \emph{ Biometrica}, vol. 38, pp. 481-482.

[6] B. V. Gnedenko and A.N.Kolomogorov (1954), \emph{Limit distributions for sums of independent random variables}, Addison-Wesely, Reading.

[7]  I. S. Gradshteyn, I. M. Ryzhik (1994), Table of integrals, series and products, $5^{th}$ Edition, Academic Press, New York.

[8] B. Hosseini (2019),Two  Metropolis-Hastings algorithms for posterior measures
with non-Gaussian priors in infinite dimensions,\emph{SIAM/ASA J. Uncertainty Quantification}, vol. 7, no 4, pp. 1185-1223.   (Also on arXiv:1804.07833 [stat.CO])

[9] M. Jeanblanc, J. Pitman and M. Yor (2002), Self-similar processes with independent increments associated with L\'evy and Bessel processes, \emph{ Stoch. Processes and their Applications} \textbf{100}, pp. 223 - 231.

[10] N. L. Johnson, S. Kotz and N. Balakrishnan (1994), Continuous univariate distributions, vol. 1, John Wiley  $\&$  Sons, Second Edition, New York.

[11] Z. J. Jurek (1996), Series of independent exponential random variables. In: \emph{Proc. 7th Japan-Russia Symposium Probab. Theory and Math. Stat.} ,  pp.174-182.
World Scientific, Singapore, New Jersey, London, Hong Kong.

[12] Z. J. Jurek (1997), Selfdecomposability: an exception or a rule?, \emph{Annales Universitatis Marie Curie-Sklodowska}, Lublin-Polonia, vol. LI. 1,  Sectio A, Mathematica,  pp. 93- 107.

[13] Z. J. Jurek (2000),  A note on gamma random variables and Dirichlet series,\emph{Stat.$\&$ Probab. Letters}, vol. 49, pp.387-392.

 [] Z. J. Jurek (2001), Remarks on the selfdecomosability and new examples, \emph{Demonstration Math.} vol. XXXIV, No 2,  pp. 241-250.

[14]  Z. J. Jurek and J.D. Mason (1993), \textit{Operator- limit Distributions
in Probability Theory}, Wiley Series in Probab. and Math.l
Statistics, New York.

[15]   Z. J. Jurek and W. Vervaat (1983), An integral representation for
selfdecomposable Banach space valued random variables,  \emph{Z. Wahrscheinlichkeits. verw. Gebiete}  \textbf{26}, 1983, pp. 247-262.
 \ \ [Also in: \ \emph{Report 8121, 25pp,  July 1981}, Katholieke Universiteit, Nijmegen, The Netherlands.]

[16] A. T. Lawrance (1982), The innovation distribution of a gamma distributed autoregressive process, \emph{Scandinavian J. Statist.}, vol. 9  (4),  pp. 234-236.

[17] M. Lo\'eve (1963), \emph{Probability theory}, Third Edition, D. Van Norstrand Company, Inc., Princeton New Jersey, Toronto, New York, London.

[18] D. N. Shanbhag, D. Pestana and M. Sreehari (1977), Some further results i infinite divisibility, \emph{Math. Proc. Camb. Phil. Soc.}, vol. 82, pp. 189-295.

% On certain self-decomposable distributions, \emph{Z. Wahrsch. verw. Gebite}, \textbf{38}, 217-222.

[19] N. G. Ushakov  (1999), \emph{Selected Topics in Characteristic Functions},  VSP, Utrecht, The Netherlands.

[20] J. G. Wendel (1961), \ The non-absolute convergence of Gil-Paleaz' inversion integral, \emph{ Ann. Math. Stat.}  vol. 32, pp. 338-339.

[21] E. T. Whittaker and  G. N. Watson (1920), A course of modern analysis, Cambridge.

\end{document}